\documentclass[12pt]{article}
\usepackage{graphicx}
\usepackage{amsmath,amsthm,amssymb,enumerate}
\usepackage{euscript,mathrsfs}
\usepackage{color}
\usepackage{dsfont}
\usepackage{url}
\usepackage[left=2cm,right=2cm,top=3.5cm,bottom=3.5cm]{geometry}
\usepackage{color}
\usepackage[framemethod=tikz]{mdframed}
\allowdisplaybreaks

\usepackage{soul}

\catcode`\@=11 \@addtoreset{equation}{section}

\catcode`\@=12

\newtheorem{Theorem}{Theorem}[section]
\newtheorem{Proposition}[Theorem]{Proposition}
\newtheorem{Lemma}[Theorem]{Lemma}
\newtheorem{Corollary}[Theorem]{Corollary}

\theoremstyle{definition}
\newtheorem{Definition}[Theorem]{Definition}

\newtheorem{Remark}[Theorem]{Remark}

\newcommand{\bTheorem}[1]{
	\begin{Theorem} \label{T#1} }
	\newcommand{\eT}{\end{Theorem}}

\newcommand{\bProposition}[1]{
	\begin{Proposition} \label{P#1}}
	\newcommand{\eP}{\end{Proposition}}

\newcommand{\bLemma}[1]{
	\begin{Lemma} \label{L#1} }
	\newcommand{\eL}{\end{Lemma}}

\newcommand{\bCorollary}[1]{
	\begin{Corollary} \label{C#1} }
	\newcommand{\eC}{\end{Corollary}}

\newcommand{\bRemark}[1]{
	\begin{Remark} \label{R#1} }
	\newcommand{\eR}{\end{Remark}}

\newcommand{\bDefinition}[1]{
	\begin{Definition} \label{D#1} }
	\newcommand{\eD}{\end{Definition}}

\newcommand{\Del}{\Delta_x}

\newcommand{\Ds}{\mathbb{D}_x}

\newcommand{\bfphi}{\boldsymbol{\varphi}}

\newcommand{\bFormula}[1]{
	\begin{equation} \label{#1}}
	\newcommand{\eF}{\end{equation}}

\newcommand{\Ov}[1]{\overline{#1}}

\newcommand{\Curl}{{\bf curl}_x}

\newcommand{\vr}{\varrho}
\newcommand{\vre}{\vr_\ep}

\newcommand{\vue}{\vu_\ep}

\newcommand{\vu}{\vc{u}}
\newcommand{\vm}{\vc{m}}

\newcommand{\vc}[1]{{\bf #1}}

\newcommand{\Div}{{\rm div}_x}
\newcommand{\Grad}{\nabla_x}

\newcommand{\dx}{\,{\rm d} {x}}

\newcommand{\dt}{\,{\rm d} t }

\newcommand{\intO}[1]{\int_{\Omega} #1 \ \dx}

\newcommand{\intTd}[1]{\int_{\mathbb{T}^d} #1 \ \dx}

\newcommand{\D}{{\rm d}}

\newcommand{\ep}{\varepsilon}

\newcommand{\Td}{\mathbb{T}^d}

\newcommand{\br}{ \nonumber \\ }

\def\softd{{\leavevmode\setbox1=\hbox{d}%
		\hbox to 1.05\wd1{d\kern-0.4ex{\char039}\hss}}}
\definecolor{Cgrey}{rgb}{0.85,0.85,0.85}
\definecolor{Cblue}{rgb}{0.50,0.85,0.85}
\definecolor{Cred}{rgb}{1,0,0}
\definecolor{fancy}{rgb}{0.10,0.85,0.10}
\definecolor{amaranth}{rgb}{0.9, 0.17, 0.31}

\newcommand\Cbox[2]{%
	\newbox\contentbox%
	\newbox\bkgdbox%
	\setbox\contentbox\hbox to \hsize{%
		\vtop{
			\kern\columnsep
			\hbox to \hsize{%
				\kern\columnsep%
				\advance\hsize by -2\columnsep%
				\setlength{\textwidth}{\hsize}%
				\vbox{
					\parskip=\baselineskip
					\parindent=0bp
					#2
				}%
				\kern\columnsep%
			}%
			\kern\columnsep%
		}%
	}%
	\setbox\bkgdbox\vbox{
		\color{#1}
		\hrule width  \wd\contentbox %
		height \ht\contentbox %
		depth  \dp\contentbox
		\color{black}
	}%
	\wd\bkgdbox=0bp%
	\vbox{\hbox to \hsize{\box\bkgdbox\box\contentbox}}%
	\vskip\baselineskip%
}

\mdfdefinestyle{MyFrame}{%
	linecolor=black,
	outerlinewidth=1pt,
	roundcorner=5pt,
	innertopmargin=\baselineskip,
	innerbottommargin=\baselineskip,
	innerrightmargin=10pt,
	innerleftmargin=10pt,
	backgroundcolor=white!20!white}

\renewcommand{\div}{{\rm div}\,}

\allowdisplaybreaks

\begin{document}


\title{\bf On physically grounded boundary conditions for the 
compressible MHD system}

\author{Jan B{\v r}ezina \and Eduard Feireisl
	\thanks{The work of E.F. was partially supported by the
		Czech Sciences Foundation (GA\v CR), Grant Agreement
		24--11034S. The Institute of Mathematics of the Academy of Sciences of
		the Czech Republic is supported by RVO:67985840.
		E.F. is a member of the Ne\v cas Center for Mathematical Modelling}
}

\date{}

\maketitle

\centerline{Faculty of Arts and Science, Kyushu University}

\centerline{744 Motooka, Nishi-ku, Fukuoka, 819-0395, Japan}

\bigskip

\centerline{Institute of Mathematics of the Academy of Sciences of the Czech Republic}

\centerline{\v Zitn\' a 25, CZ-115 67 Praha 1, Czech Republic}

\begin{abstract}
	
We consider a general compressible MHD system, where the magnetic field propagates in a heterogeneous medium. Using suitable penalization in terms of the transport coefficients we perform several singular limits. As a result we obtain:
\begin{enumerate}
	\item A rigorous justification of physically grounded boundary conditions for the compressible MHD system on a bounded domain.
	\item Existence of weak solutions for arbitrary finite energy initial data in the situation the Maxwell induction equation holds also outside the fluid domain.
	\item A suitable theoretical platform for numerical experiments on domains with geometrically 	complicated boundaries.
\end{enumerate}	

\end{abstract}


{\small

\noindent
{\bf 2020 Mathematics Subject Classification:}{ 
(primary); 
(secondary) }

\medbreak
\noindent {\bf Keywords:}

\tableofcontents

}

\section{Introduction}
\label{i}

The magnetohydrodynamics (MHD) system of equations couples {fluid equations} with a simplified Maxwell system describing the propagation of {a} magnetic field. The fluid {motion is driven} by 
the Lorentz force that can be conveniently written as a divergence of 
Maxwell's tensor, while the electric current in Ohm's law in the induction equation is augmented by the vector product of the fluid velocity with the magnetic {field} vector. If the fluid in question 
is compressible barotropic, the relevant system of equations reads 
	\begin{align}
	\partial_t \vr + \Div (\vr \vu) &= 0,\ \label{i1} \\
	\partial_t (\vr \vu) + \Div (\vr \vu \otimes \vu) + 
	\Grad p &= \Div \mathbb{S} - \beta \vu - \Div \left[ \mu \left( \vc{H} \otimes \vc{H} - \frac{1}{2} |\vc{H}|^2 \mathbb{I} \right) \right] \label{i2}, \\
	\partial_t (\mu \vc{H}) + \Curl ( \mu \vc{H} \times \vu) + 
	\Curl (\eta \Curl \vc{H}) &= 0,\ \Div (\mu \vc{H}) = 0, \label{i3}
\end{align}  
with the unknowns $\vr = \vr(t,x)$ - the fluid mass density, 
$\vu = \vu(t,x)$ the velocity of the fluid, and $\vc{H} = \vc{H}(t,x)$ the magnetic {field} vector. The pressure $p = p(\vr)$ is given constitutively, the viscous stress $\mathbb{S} = \mathbb{S}(\Ds \vu)$ obeys Newton's rheological law, 
\begin{equation} \label{i4}
	\mathbb{S}(\Ds \vu) = \nu \left( \Grad \vu + \Grad^t \vu - \frac{2}{d} 
	\Div \vu \mathbb{I} \right) + \lambda \Div \vu \mathbb{I},\ 
	\Ds \vu \equiv  \frac{1}{2} \left( \nabla \vu + \nabla^t \vu \right).
\end{equation}
In addition, the equations \eqref{i2}, \eqref{i3} depend on transport coefficients:

$\beta$ \dotfill friction parameter 

$\nu$ \dotfill shear viscosity coefficient 

$\lambda$ \dotfill bulk viscosity coefficient

$\mu$ \dotfill magnetic permeability 

$\eta$ \dotfill electric resistivity 

\noindent They are given non--negative  functions of the spatial 
variable $x$. In addition, we suppose $\nu$, $\mu$, and $\eta$ are strictly positive bounded below away from zero. 

While the fluid is usually confined to a well specified bounded domain with an impermeable boundary, the induction equation \eqref{i3} 
is in principle satisfied outside the fluid domain with $\vu = 0$, and with possibly different {values of} coefficients $\mu$ and $\eta$. The relevant boundary conditions specifying the value of the field $\vc{H}$ on the boundary of the fluid domain are therefore of \emph{transmission} type.
This approach was adopted in the seminal work of Ladyzhenskaya and Solonnikov \cite{LadSol} in the context of incompressible MHD. Later,  
Meir and Schmidt \cite{MeiSc1}, \cite{MeiSc2} developed a similar method 
applicable to stationary incompressible MHD problems. To the best of our knowledge, this approach has never been applied to the compressible MHD models. One of the principal objectives of this paper is to fill this gap.

Solving the induction equation outside the fluid domain may be technically complicated, in particular in view of numerical simulations. For this reason, the system \eqref{i1}--\eqref{i3} is often supplemented by explicit boundary conditions for the field $\vc{H}$. 
Typically, the so--called perfect electric conductor (PEC) boundary conditions are the most frequently used choice.
As pointed out by Shercliff \cite{Sher}:
\begin{quote}
{\it A perfectly conducting solid boundary is so much favoured by theoreticians because it shields  the fluid from external conditions and eliminates consideration of all but the fluid region. This selection usually constitutes escapism, however,  justifiable only as a simplifying first approximation.}	
	
\end{quote}	

What is more, a proper choice of the boundary conditions for $\vc{H}$ is a delicate issue. Assuming $\mu$ and $\eta$ constant we may rewrite the differential operator in {\eqref{i3}} as 
\[
\Curl \Curl \vc{H} = \Grad \Div \vc{H} - \Del \vc{H} = - \Del \vc{H}.
\]
One is therefore tempted to adopt the homogeneous Dirichlet boundary 
condition $\vc{H}|_{\partial \Omega} = 0$. Unfortunately, as observed by 
Lassner \cite{Lass} (see also Str\" ohmer \cite{Stroh}), the problem 
{\eqref{i3}} becomes overdetermined unless a variant of ``pressure'' is added in the formulation of Ohm's law. This fact is often incorrectly interpreted as a discrepancy between the weak and strong solutions of the problem. {Therefore our} second objective is to justify most of the commonly used boundary conditions as singular limits of the transmission conditions for a suitable scaling of {the} transport coefficients.

For the sake of simplicity and also in view of possible numerical computations, we start with the compressible MHD system with the space periodic boundary conditions, meaning the spatial domain can be identified with a flat torus 
\[
\Td = \left( [-L,L] |_{\{ -L, L\} } \right)^d,\ d=2,3, 
\]
where $L$ is large enough to absorb the fluid domain. At this stage, we consider the transport coefficients to be continuously differentiable functions of $x \in \Td$.

Motivated by the seminal work of Ladyzhenskaya and Solonnikov \cite{LadSol}, we fix a Lipschitz domain 
\begin{equation} \label{i5}
	\Omega \subset R^d ,\ \Omega \subset \Ov{\Omega} \subset 
	(-L,L)^d.
\end{equation} 
In addition, we consider another Lipschitz domain 
\begin{equation} \label{i6}
	\Omega_{\rm int} \subset \Ov{\Omega}_{\rm int} \subset \Omega, 
\end{equation}
representing a solid object immersed in the fluid. Accordingly, the fluid occupies the domain 
\begin{equation} \label{i6a}
{\Omega_F = \Omega \setminus \Ov{\Omega}_{\rm int}.} 
\end{equation}
Finally, we set 
\begin{equation} \label{i6b}
{\Omega_{\rm ext} = \Td \setminus \Ov{\Omega}.}
\end{equation}

 \begin{figure}[h]
 \begin{center}
 \includegraphics[scale=0.8]{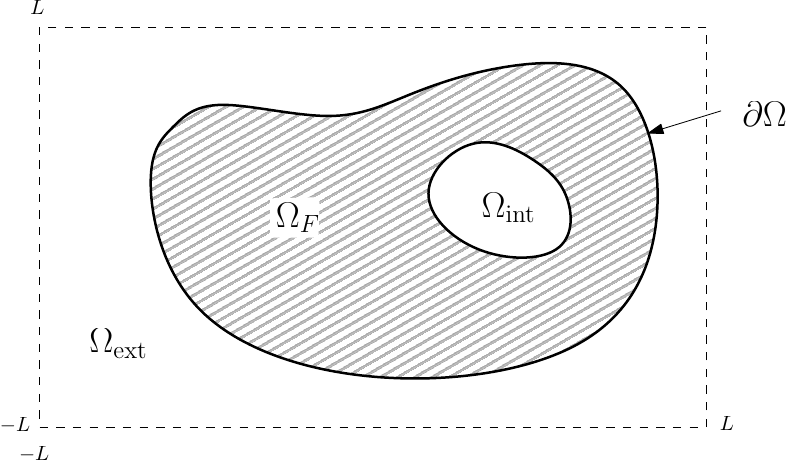}
\end{center}
 \end{figure}

Our goal is to show {the} existence of weak solutions to several initial--boundary value {problems} for the compressible MHD system resulting from appropriate penalization of the transport coefficient. 
Specifically, we consider the following problems:
\begin{itemize}
	\item {The equation of continuity \eqref{i1} holds in the whole domain $(0,T) \times \Td$}.
	\item The fluid {momentum equation \eqref{i2} is} satisfied in $(0,T) \times \Omega_F$, 
with $\beta = 0$ and constant viscosity coefficients $\nu = \nu_F > 0$, $\lambda = \lambda_F \geq 0$, supplemented with the no--slip boundary condition 
\begin{equation} \label{i7}
	\vu|_{\partial \Omega_F} = 0.
\end{equation}	

\item The induction equation {\eqref{i3}} is satisfied in $(0,T) \times \Omega$, with constant but generally different transport coefficients 
\[
\mu = \left\{ \begin{array}{l} \mu_{\rm int} > 0 \ \mbox{in}\ \Omega_{\rm int} \\ {\mu_F > 0} \ \mbox{in}\ \Omega_F \end{array} \right., 
\eta = \left\{ \begin{array}{l} \eta_{\rm int} > 0 \ \mbox{in}\ \Omega_{\rm int} \\ {\eta_F > 0} \ \mbox{in}\ \Omega_F \end{array} \right.,
\]
{where the fluid velocity $\vu$ is set to be zero outside $\Omega_F$},
with the corresponding transmission boundary conditions on $\partial \Omega_{\rm int}$.
\item {The magnetic field vector} $\vc{H}$ will satisfy one of the following boundary conditions on $\partial \Omega$:
\begin{enumerate}
	\item {\bf Perfect electric isolator.}
	\begin{equation} \label{i8}
		\vc{H} \times \vc{n} = \vc{H}_{\rm ext} \times \vc{n},\ 
		\mu_F \vc{H} \cdot \vc{n} = \mu_{\rm ext} \vc{H}_{\rm ext} \cdot \vc{n}
		\ \mbox{in}\ (0,T) \times \partial \Omega,
	\end{equation}	
where $\vc{H}_{\rm ext}$ is an external magnetic field 
satisfying \begin{align}
	\Curl \vc{H}_{\rm ext} = 0 ,\ \Div (\mu_{\rm ext} \vc{H}_{\rm ext} ) = 0	\ \mbox{in}\ (0,T) \times \Omega_{\rm ext}
	\nonumber
\end{align}
	
These are the ``non--local'' boundary conditions considered in problem II by Ladyzhenskaya and Solonnikov \cite{LadSol}, see also Sakhaev and Solonnikov \cite{SakSol}.

\item {\bf Perfect magnetic conductor (PMC).}
\begin{equation} \label{i9}
	\vc{H} \times \vc{n} = 0 \ \mbox{in}\ (0,T) \times \partial \Omega
\end{equation}	
{Note that these are the \emph{correct} Dirichlet type boundary conditions to be imposed on $\vc{H}$. Here, $\Omega_{\rm ext}$ is a  perfect magnetic conductor with $\mu_{\rm ext} \to \infty$.}

\item
{\bf Perfect electric conductor (PEC).}

\begin{equation} \label{i10}
	\vc{H} \cdot \vc{n} = 0,\  
	\Curl \vc{H} \times \vc{n} = 0 \ \mbox{in}\ (0,T) \times \partial \Omega
\end{equation}

These {are} probably the most frequently used boundary conditions in mathematical treatments of MHD. {The external domain $\Omega_{\rm ext}$} is a perfect conductor (possibly superconductor), for which	 
 \[
	\mu_{\rm ext} \to 0,\ \eta_{\rm ext} \to 0.	
\]

\item {\bf Isolator type boundary conditions.}

Str\" ohmer \cite{Stroh} proposed the boundary conditions
\begin{equation} \label{i11}
	\vc{H} \cdot \vc{n} = 0,\  
	\Curl \vc{H} \cdot \vc{n} = 0 \ \mbox{in}\ (0,T) \times \partial \Omega.	
\end{equation}
The inhomogeneous version of {\eqref{i11}} has been recently examined by Neustupa, Perisetti, and Yang \cite{NePeYa}. 

Although mathematically {still} admissible, these {boundary conditions} correspond to a rather unrealistic material $\Omega_{\rm ext}$ satisfying 
\[
	\mu_{\rm ext} \to 0, \ \eta_{\rm ext} \to \infty.	  
\]

\end{enumerate}

\end{itemize}	

The paper is organized as follows. In Section \ref{P}, we introduce the starting primitive system defined on the periodic spatial domain {$\Td$}, with smooth transport coefficients. We recall the concept of weak solution both for the primitive system and the asymptotic singular limit problems.
{In} Section \ref{Z}, we consider {various} types of penalizations of the transport coefficients. Having collected the necessary preliminary material we state our main results in Section \ref{M}. Section \ref{A} 
is devoted to the proof of {the} main results.

\section{Problem with smooth transport coefficients}
\label{P}

We start with the compressible MHD system with spatially periodic boundary conditions and spatially dependent transport coefficients
defined on the flat torus
\[
\Td = \left( [-L,L] |_{\{ -L, L\} } \right)^d. 
\]
For the sake of simplicity, we consider the isentropic pressure state equation
\[
p(\vr) = a \varrho^\gamma,\ a > 0, \ \gamma > \frac{d}{2}. 
\]
More complicated equations of state including a non--monotone pressure can be handled by essentially the same approach, see e.g. \cite[Chapter 6]{EF70}.

\subsection{Weak formulation of the primitive system}

{We suppose} the transport coefficients belong to the class $C^1(\Td)$, and, in addition, 
\begin{equation} \label{P1}
	\beta \geq 0,\ \nu > 0, \lambda \geq 0,\ \mu > 0,\ \eta > 0 
	\ \mbox{in}\ \Td.
\end{equation}

\begin{Definition}[\bf Weak solutions to the primitive system]\label{DP1}
We shall say that a trio of functions $(\vr, \vu, \vc{H})$ is \emph{weak solution} of the compressible MHD system \eqref{i1}--\eqref{i3} in the space--time cylinder $(0,T) \times \Td$, with the initial data $(\vr_0, \vm_0, \vc{H}_0)$ if the following holds:	
\begin{itemize}
\item {\bf Equation of continuity.}
The functions $\vr$, $\vu$ belong to the class $\vr \geq 0$, {$\vr \in C_{\rm weak}\left( 
[0,T]; L^\gamma(\Td) \right) $, $\vu \in L^2(0,T; W^{1,2}(\Td; R^d))$,} and the integral identity
\begin{equation} \label{P2} 
\intTd{ \vr (\tau, \cdot) \varphi } = \intTd{ \vr_0 \varphi } +  \int_0^\tau \intTd{ \vr \vu \cdot \Grad \varphi } \dt,\ 
0 \leq \tau \leq T
\end{equation}
holds for any $\varphi \in C^1(\Td)$. In addition, the renormalized equation
\begin{equation} \label{P3} 
	\intTd{ b(\vr) (\tau, \cdot) \varphi } = \intTd{ b(\vr_0) \varphi } +  \int_0^\tau \intTd{ \left[ b(\vr) \vu \cdot \Grad \varphi + 
		\Big( b(\vr) - b'(\vr) \vr \Big) \Div \vu \varphi \right]} \dt
\end{equation}
holds for any $\varphi \in C^1(\Td)$, and any $b \in C^1[0, \infty)$,
$b' \in C_c[0, \infty)$.

\item {\bf Momentum equation.} 
The momentum $\vr \vu \in C_{\rm weak}([0,T]; L^{\frac{2 \gamma}{\gamma + 1}}(\Td; R^d))$, the magnetic {field} vector $\vc{H} \in L^\infty(0,T; L^2(\Td; R^d))$, and the integral identity 
\begin{align}
	\intTd{\vr \vu (\tau, \cdot) \cdot \bfphi } &= \intTd{\vm_0 \bfphi} + \int_0^\tau \intTd{ \Big[
		(\vr \vu \otimes \vu): \Grad \bfphi + p (\vr) \Div \bfphi \Big] } \dt \br
	&- \int_0^\tau \intTd{ \mathbb{S} (\Ds \vu) : \Grad \bfphi } \dt +
	\int_0^T \intTd{ \beta \vu \cdot \bfphi } \dt \br &+
	\int_0^\tau \intTd{ \mu \Big[ \vc{H} \otimes \vc{H} - \frac{1}{2} |\vc{H}|^2 \mathbb{I} \Big] : \Grad \bfphi } \dt
	\label{P4} 		
\end{align}	
{holds} for any $0 \leq \tau \leq T$ and any $\bfphi \in C^1(\Td; R^d)$.
\item {\bf Induction equation.} {The magnetic field vector 
$\vc{H} \in C_{\rm weak}([0,T]; L^2(\Td; R^d))$, $\Curl \vc{H} 
\in L^2(0,T; L^2(\Td; R^d))$}, and the integral identity
\begin{equation} \label{P5}
	\intTd{ \mu \vc{H}(\tau, \cdot) \cdot \bfphi} = 
	\intTd{\mu \vc{H}_0 \cdot \bfphi }
	-
	\int_0^\tau \intTd{ \Big[  
		(\mu \vc{H} \times \vu) \cdot \Curl \bfphi + \eta \Curl \vc{H} \cdot \Curl \bfphi \Big] } \dt 
\end{equation}		
{holds} for any $0 \leq \tau \leq T$, and any $\bfphi \in C^1(\Td)$. Moreover, we suppose 
\begin{equation} \label{P6}
\intTd{ \mu \vc{H}_0 \cdot \Grad \phi } = 0 \ \mbox{for all}\ \phi \in C^1(\Td) 
\end{equation}
yielding 
\begin{equation} \label{P7} 
	\intTd{ \mu \vc{H}(\tau, \cdot) \cdot \Grad \phi } = 
0 \ {\mbox{meaning}}\  \Div (\mu \vc{H}(\tau, \cdot)) = 0 ,\ 
 {\mbox{in the sense of distributions.}}
\end{equation}	

\end{itemize}	
	\end{Definition}
	
\begin{Remark} \label{RP1}
	
We recall Gaffney inequality, see e.g. {Csato, Dacorogna, Sil \cite{CsDaSi}}, 
\begin{equation} \label{P8}	 
\| \Grad \vc{h} \|_{L^2(Q; R^{d \times d})}	\leq c(Q) 
\Big[ \| \vc{h} \|_{L^2(Q; R^d)} + 
\| \Curl \vc{h} \|_{L^2(Q; R^d)} + \| \Div \vc{h} \|_{L^2(Q)} \Big]
\end{equation}	
for any $\vc{h} \in W^{1,2}(Q; R^d)$ satisfying either $\vc{h} \times \vc{n}|_{\partial Q} = 0$ or $\vc{h} \cdot \vc{n}|_{\partial Q} = 0$. It is easy to see {that \eqref{P8} remains valid 
on $\Td$}. Consequently, as $\mu$ is smooth {and strictly positive}, and $\vc{H}$ satisfies \eqref{P7}, we may infer that any weak solution in the sense of Definition \ref{DP1} belongs to the class
\begin{equation} \label{P9}
	\vc{H} \in L^2(0,T; W^{1,2}(\Td; R^d)).
\end{equation}

	\end{Remark}

\subsection{Energy balance}

We focus on the weak solutions satisfying some form of energy balance. 
Let $P$ be the so--called pressure potential satisfying 
\[
P'(\vr)\vr - P(\vr) = p(\vr). 
\]
It is convenient to introduce the momentum $\vm = \vr \vu$ to write the total energy of the system in the form 
\begin{equation} \label{P10}
	E(\vr, \vm, \vc{H}) = \left\{ \begin{array}{l} \frac{1}{2} \frac{|\vm|^2}{\vr} + P(\vr) + \frac{1}{2} \mu |\vc{H}|^2 \ \mbox{if}\ 
		\vr > 0 \\ \\
		P(0) + \frac{1}{2} \mu |\vc{H}|^2 \ \mbox{if}\ \vr = 0, \vm = 0, \\ \\
		\infty \ \mbox{otherwise.}	 
	\end{array} \right.		
\end{equation}	
If the pressure $p$ is an increasing function of the density, the energy 
is a convex l.s.c. function defined on $R^{2d +1}$.

The {primitive} system is augmented by the energy inequality 
\[ 
\frac{\D }{\dt} \intTd{ E(\vr ,\vm, \vc{H}) } + 
	\intTd{ \left[ \mathbb{S}(\Ds \vu) : \Ds \vu + \eta |\Curl \vc{H}|^2 \right] } + \intTd{ \beta |\vu|^2 } \leq 0
	\]
that can be written in the weak form as 
	\begin{align} 
		- \int_0^T \partial_t \psi \intTd{ E(\vr, \vm , \vc{H}) } \dt &+ 
		\int_0^T \psi \intTd{ \left[ \mathbb{S}(\Ds \vu) : \Ds \vu + \eta |\Curl \vc{H}|^2 \right] } \dt \br &+ \int_0^T \psi \intTd{ \beta |\vu|^2 } \dt  
		\leq \psi(0) \intTd{ E(\vr_0, \vm_0 , \vc{H}_0) }
		\label{P11}
	\end{align}	
	for any $\psi \in C^1_c[0,T)$, $\psi \geq 0$.

Relation \eqref{P11} can be integrated in time to obtain its weaker form

\begin{align} \label{P12}
	\intTd{ E(\vr, \vm , \vc{H}) (\tau, \cdot) } &+ 
	\int_0^\tau \intTd{ \left[ \mathbb{S}(\Ds \vu) : \Ds \vu + \eta |\Curl \vc{H}|^2 \right] } \dt \br &+ \int_0^\tau  \intTd{ \beta |\vu|^2 } \dt  
	\leq \intTd{ E(\vr_0, \vm_0 , \vc{H}_0) },\ 0 \leq \tau \leq T.
\end{align}

\begin{Definition}[\bf Energy weak solution] \label{DP2}
\begin{itemize} 
	\item
We say that a weak solution $(\vr, \vu, \vc{H})$ of the compressible MHD system {\eqref{i1}--\eqref{i3}} is \emph{finite energy weak solution} if it satisfies the energy inequality \eqref{P12}.
\item We say that a weak solution $(\vr, \vu, \vc{H})$ of the compressible MHD system {\eqref{i1}--\eqref{i3}} is \emph{dissipative weak solution} if it satisfies  the energy inequality {in the differential form} \eqref{P11}.
\end{itemize}
\end{Definition}

\subsection{Existence and regularity of the weak solutions}

If the initial data belong to the class 
\[
\intTd{ E(\vr_0, \vm_0, \vc{H}_0) } < \infty, 
\intTd{ \vr_0 } = M > 0,\ \Div (\mu \vc{H}_0) = 0,
\]
it is possible to show {the existence} of global in time, meaning defined on a given time interval $(0,T)$, dissipative weak solutions of the compressible MHD system in the space--time cylinder $(0,T) \times \Td$. The approach developed in \cite{DUFE2} is easy to adapt to the case of 
smooth transport coefficients. Note that the transport coefficients in 
\cite{DUFE2} depend on the temperature which {is more difficult} than the present case. 

Regularity (integrability) of weak solutions is essentially given by the energy estimates \eqref{P12}. As $\mu$ is positive {and} continuously differentiable, we can use {the} argument of Remark \ref{RP1} to conclude
\[
\vc{H} \in L^\infty(0,T; L^2(\Td; R^d)) \cap 
L^2(0,T; W^{1,2}(\Td; R^d)),
\]
cf. also \cite{DUFE2}.

\section{Penalization of transport coefficients}
\label{Z}

The penalization of the transport coefficients may be put in a unified framework as follows:
For each transport coefficient 
\[
\Lambda = \beta, \nu, \lambda, \mu, \eta, 
\]
we prescribe three constant values 
\[
\Lambda^\ep_F,\ \Lambda^\ep_{\rm int}, \ \Lambda^\ep_{\rm ext}
\]
together with a family of smooth approximations {enjoying the following properties:}
\begin{align}
	\Lambda_\ep &\in C^1(\Td),\ \min\{ \Lambda^\ep_F, 
	\Lambda^\ep_{\rm int}, \Lambda^\ep_{\rm ext} \} \leq \Lambda_\ep  \leq \max\{ \Lambda^\ep_F, 
	\Lambda^\ep_{\rm int}, \Lambda^\ep_{\rm ext} \} , \br	
	\Lambda_\ep(x) &= \left\{ \begin{array}{l} \Lambda^\ep_{\rm int} \ \mbox{if}\ x \in {\Omega_{\rm int}}\\ \Lambda_{\rm ext}^{\ep} \ \mbox{if}\ x \in {\Omega_{\rm ext}}  \end{array} \right. , \br 
	\mbox{for any compact}&\ Q \subset \Omega_F,\ \mbox{there exists}\ \ep_0 
	\ \mbox{such that}\ \Lambda_\ep(x) = \Lambda^\ep_F  \ \mbox{if}\ x \in Q
	\ \mbox{and}\ 0 < \ep < \ep_0.
	\label{Z1}	
\end{align}

\subsection{Singular limit in the fluid equations}

Our approach is based on a proper choice of scaling of the parameters 
$(\Lambda^\ep_F, \Lambda^\ep_{\rm int}, \Lambda^\ep_{\rm ext})$ and performing the limit $\ep \to 0$. In this process, we \emph{always} assume the initial energy 
\begin{equation} \label{Z2}
\intTd{ E(\vr_{0,\ep}, \vm_{0, \ep}, \vc{H}_{0, \ep}) } \leq c, 
\end{equation}
where the constant {$c$} is independent of $\ep$. In such a way, the energy inequality \eqref{P12} yields uniform bounds independent of $\ep$.

The first goal is to keep the fluid confined to the domain $\Omega_F$. To this end, we consider the standard Brinkmann type penalization 
\begin{equation} \label{Z3}
	\beta^\ep_F = 0,\ \beta^\ep_{\rm int} \to \infty, \
	\beta^\ep_{\rm ext} \to \infty.
\end{equation}
In addition, we consider the high viscosity limit outside $\Omega_F$, meaning 
\begin{align} 
\nu^\ep_F &= \nu_F > 0,\ \lambda^\ep_F = \lambda_F \geq 0,\br 
\nu^\ep_{\rm int} &\to \infty,\ \nu^\ep_{\rm ext} \to \infty \ \mbox{if}\ d = 3, \br 	
\lambda^\ep_{\rm int} &\to \infty,\ \lambda^\ep_{\rm ext} \to \infty 
\ \mbox{if}\ d = 2.
\label{Z4}
\end{align}
The remaining viscosity coefficients in the cases $d = 2,3$ can be arbitrary as long as they remain non--negative. The necessity of having large bulk viscosity in the 2-d case is due to the gap in the corresponding Korn inequality, cf. Dain \cite{DAIN}.

\subsubsection{Singular limit in the equation of continuity}

If both \eqref{Z3} and \eqref{Z4} are satisfied, we can use the traceless variant of Korn's inequality (cf. \cite[Chapter 11, Section 10]{FeNo6A}) together with the energy bounds \eqref{P12} to conclude 
\begin{align} 
\vue &\to \vu \ \mbox{weakly in}\ L^2(0,T; W^{1,2}(\Td; R^d)), \br
\vue &\to 0 \ \mbox{(strongly) in}\ {L^2(0,T; L^2(\Omega_{\rm int} \cup \Omega_{\rm ext}; R^d)}, \br
\Div \vue &\to 0 \ \mbox{in}\ {L^2(0,T; L^2(\Omega_{\rm int} \cup \Omega_{\rm ext}))}
\label{Z5}	
	\end{align}
passing to a suitable subsequence as the case may be. In addition, we may use the energy bounds together with the equation of continuity \eqref{P2} to deduce 
\begin{equation} \label{Z6}
	\vre \to \vr \ \mbox{in}\ C_{\rm weak}([0,T]; L^{\gamma}(\Td)) 
\end{equation}
again for a suitable subsequence. The limit satisfies the equation of continuity \eqref{P2} in $(0,T) \times \Td$. As $\vu = 0$ in $\Omega_{\rm int} \cup \Omega_{\rm ext}$ we have $\vr = \vr_0$ in 
{$\Omega_{\rm int} \cup \Omega_{\rm ext}$} if the initial density is independent of $\ep$.

Finally, we may use the renormalized equation of continuity \eqref{P3} 
and the strong convergence of {$\vue$ and $\Div \vue$} on the solid part established in \eqref{Z5} to conclude 
\[
\intTd{ b(\vre) \varphi } \to \intTd{ b(\vr_0) \varphi } 
\ \mbox{uniformly in}\ C[0,T] 
\]
for any {$\varphi \in C^1_c(\Omega_{\rm int} \cup \Omega_{\rm ext})$}, which yields strong convergence of the densities on the solid part, 
\begin{equation} \label{Z7}
\vre \to {\vr} \ \mbox{in}\ {C([0,T]; L^1(\Omega_{\rm int} \cup \Omega_{\rm ext}))}.	
\end{equation}

\subsubsection{Singular limit in the momentum equation}

In order to perform the limit in the momentum equation, we have to specify the transport coefficients $\mu_\ep$, $\eta_\ep$ at least in the domain $\Omega$. We set 
\begin{equation} \label{Z8}
\mu^\ep_{\rm int} = \mu_{\rm int} > 0, 
\eta^\ep_{\rm int} = \eta_{\rm int} > 0,\ \mu^\ep_F = \mu_F > 0,\ 
{\eta^\ep_F = \eta_F > 0},
\end{equation}	
leaving the specific choice of $\mu^\ep_{\rm ext}$, $\eta^\ep_{\rm ext}$ open. As the test functions in the {momentum equation} are compactly supported {in $\Omega_F$}, the piece of information \eqref{Z8} is {sufficient to establish the limit in the momentum equation in $\Omega_F$.}

As the initial energy is bounded, we deduce from \eqref{P12} 
\begin{equation} \label{Z9}
	\vc{H}_\ep \to \vc{H}\ \mbox{weakly-(*) in}\ L^\infty(0,T; L^2(Q; R^d)) \ \mbox{for any compact} \ Q \subset \Omega.
\end{equation}	
In addition, 
\begin{equation} \label{Z10}
	\Curl \vc{H}_\ep \to \Curl \vc{H} \ \mbox{weakly in}\ 
L^2(0,T; L^2(Q; R^d)) \ \mbox{for any compact}\ Q \subset \Omega 
\end{equation}
at least for a suitable subsequence. Moreover, the {limit belongs to the space} 
\begin{equation} \label{Z11}
	\vc{H} \in L^\infty(0,T; L^2(\Omega; R^d)),\ 
	\Curl \vc{H} \in L^2(0,T; L^2(\Omega; R^d)).	 
\end{equation}

Going back to the induction equation, the uniform bounds established above further imply
\begin{align} 
\mu^\ep \vc{H}_\ep &\to \mu_{\rm int} \vc{H} \ \mbox{in} \ C_{\rm weak}([0,T]; L^2(\Omega_{\rm int}; R^d)), \br 
\mu^\ep \vc{H}_\ep &\to \mu_F \vc{H} \ \mbox{in} \ C_{\rm weak}([0,T]; L^2(Q; R^d)) \ \mbox{for any compact}\ Q \subset \Omega_F.  
\label{Z12}
\end{align}	
This yields
\begin{equation} \label{Z13}
\vc{H}_\ep \to \vc{H} \ \mbox{in}\ C_{\rm weak}([0,T]; L^2(Q; R^d)) \ \mbox{for any compact}\ Q \subset {\Omega}.
\end{equation}	
Finally, by virtue of Gaffney inequality \eqref{P8}, 
\begin{align} 
\vc{H}_\ep &\to \vc{H} \ \mbox{weakly in}\ L^2(0,T; W^{1,2}(Q; R^d)) \ \mbox{for any compact}\ Q \subset \Omega_F , \br
\vc{H}_\ep \otimes \vc{H}_\ep - \frac{1}{2}|\vc{H}_\ep|^2 \mathbb{I} &\to \vc{H} \otimes \vc{H} - \frac{1}{2}|\vc{H}|^2 \mathbb{I}
\ \mbox{weakly in}\ L^q((0,T) \times Q; R^{d \times d})\ \mbox{for any compact}\ Q \subset \Omega_F , \br
\mbox{and some}\ q &> 1.
	\label{Z14}
\end{align}

With the strong convergence of Maxwell's stress established in \eqref{Z14}, the limit passage in the equation of momentum in $\Omega_F$ is a routine matter, 
see e.g. \cite{BaFeLMMiYu} or \cite{FeNeSt}. 
We recall the main steps for reader's convenience:

\begin{enumerate}
	\item {\bf Pressure estimates.} The available energy estimates yield the pressure sequence $(p(\vre))_{\ep > 0}$ bounded in a non--reflexive space {$L^\infty(0,T; L^1(\Omega_F))$}. In order to improve integrability of the pressure in space, one may use 
	\[
	\bfphi = \phi \Grad \Delta^{-1} [ \vr^\alpha ] 
	\ \mbox{for some}\ \alpha > 0,\ \phi \in C^\infty_c(\Omega_F)
	\]
	as a test function in the momentum equation {\eqref{P4}}. After a bit tedious but nowadays well understood series of estimates, see e.g. 
	\cite[Chapter 5]{EF70}, we obtain uniform local estimates 
	\begin{equation} \label{aA1}
		\| p(\vre) \vre^\alpha \|_{L^1((0,T) \times Q)} \leq c(Q) 
		\ \mbox{for any compact}\ Q \subset \Omega_F.	
	\end{equation}
	Here, the crucial observation is that {Maxwell's tensor satisfies \eqref{Z14}}, in particular, it is $q-$integrable on any compact $Q$.
	
	\item {\bf Weak limit in the {momentum} equation.}
	Using \eqref{Z14} and the pressure estimates claimed in the preceding step, we may pass to the limit in the momentum equation {\eqref{P4}} obtaining 
	\begin{align}
		\int_0^T & \int_{\Omega_F} \Big[ \vr \vu \cdot \partial_t \bfphi + 
		(\vr \vu \otimes \vu): \Grad \bfphi + \Ov{p(\vr)} \Div \bfphi \Big]  \dx \dt \br
		&= \int_0^T \int_{\Omega_F} {\mathbb{S}(\Ds \vu)} : \Grad \bfphi \dx \dt \- 
		\int_0^T \int_{\Omega_F} \mu_F \Big[ \vc{H} \otimes \vc{H} - \frac{1}{2} |\vc{H}|^2 \mathbb{I} \Big] : \Grad \bfphi \dx \dt \br &- \int_{\Omega_F} \vr_0 \vu_0 
		\cdot \bfphi(0, \cdot) \dx 
		\label{aA2} 		
	\end{align}	
	for any $\bfphi \in C^1_c([0, T) \times \Omega_F; R^d)$. Recall 
	{that test functions are compactly supported in $\Omega_F$}. Here, the symbol $\Ov{p(\vr)}$ stands for the weak limit of $(p(\vre))_{\ep > 0}$. 
	
	\item {\bf Strong convergence of densities}
	The final step consists in showing 
	\[ 
	\Ov{p(\vr)} = p(\vr). 
	\]
	Note that once this is established, we also have {$p(\vr) \in L^\infty((0,T); L^1(\Omega_F))$.}	
	This is, of course, the most delicate argument based on Lions' theory 
	\cite{LI4} of weak solutions of the compressible Navier--Stokes system and its extension in \cite{EF70}. The key point is showing Lions' identity 	 
	\begin{equation} \label{aA3}
		{ \Ov{ \Big[ p(\vr) + (\lambda + 2 \nu) \div \vu \Big] b(\vr) } = 
		\Ov{ \Big[ p(\vr) + (\lambda + 2 \nu) \div \vu \Big]} \ \Ov{b(\vr)} } 
	\end{equation}
	for any $b \in C^1[0, \infty)$, $b' \in C_c([0, \infty))$ in 
	$(0,T) \times Q$ for any compact {$Q \subset \Omega_F$}. From this, one deduces the oscillations defect measure bound 	
	\begin{equation} \label{aA4}
		{\rm osc}_{\gamma + 1}((0,T) \times \Omega_F)[\vre \to \vr] 
		= \sup_{k \geq 1} \Big( \limsup_{\ep \to 0} \| T_k (\vre) - T_k(\vr) \|_{L^{\gamma + 1}((0,T) \times \Omega_F)} \Big) < \infty. 
	\end{equation}
	Note carefully that \eqref{aA4} holds on the whole fluid domain $\Omega_F$ while the relation \eqref{aA3} is local, see \cite{FNP1}. 
	With \eqref{aA4} at hand, the proof of strong convergence of the densities can be performed exactly as in \cite[Chapter 6]{EF70}.

\end{enumerate}

Let us summarize the results concerning convergence in the fluid domain. 
\begin{Proposition} \label{ZP1}
	Under the scaling \eqref{Z3}, \eqref{Z4}, \eqref{Z8}, we have	
\begin{align}	
\vue &\to \vu \ \mbox{weakly in}\ L^2(0,T; W^{1,2}(\Td; R^d)), \br
\vue &\to 0 \ \mbox{(strongly) in}\ L^2(0,T; L^2(\Omega_{\rm int} \cup {\Omega_{\rm ext}}; R^d), \br
\Div \vue &\to 0 \ \mbox{in}\ L^2(0,T; L^2(\Omega_{\rm int} \cup {\Omega_{\rm ext}})), \br
	\vre &\to \vr \ \mbox{in}\ C_{\rm weak}([0,T]; L^{\gamma}(\Td)), \br 
\vre &\to \vr_0 \ \mbox{in}\ C([0,T]; L^1(\Omega_{\rm int} \cup {\Omega_{\rm ext}})), \br 
	\vc{H}_\ep &\to \vc{H}\ \mbox{weakly-(*) in}\ L^\infty(0,T; L^2(Q; R^d)) \ \mbox{for any compact} \ Q \subset \Omega, {\mbox{cf. also \eqref{Z13}}}, \br
		\Curl \vc{H}_\ep &\to \Curl \vc{H} \ \mbox{weakly in}\ 
	L^2(0,T; L^2(Q; R^d)) \ \mbox{for any compact}\ Q \subset \Omega
\label{Z19}
\end{align}
passing to a suitable subsequence as the case may be. The limit satisfies the equation of continuity
\begin{equation} \label{Z20} 
	\intTd{ \vr (\tau, \cdot) \varphi } = \intTd{ \vr_0 \varphi } +  \int_0^\tau \intTd{ \vr \vu \cdot \Grad \varphi } \dt,\ 
	0 \leq \tau \leq T
\end{equation}
for any $\varphi \in C^1(\Td)$, and the momentum equation
\begin{align}
	\int_{\Omega_F} \vr \vu (\tau, \cdot) \cdot \bfphi \dx &= \int_{\Omega_F} \vm_0 \cdot \bfphi \ \dx + \int_0^\tau \int_{\Omega_F} \Big[
		(\vr \vu \otimes \vu): \Grad \bfphi + p (\vr) \Div \bfphi \Big] \dx \dt \br
	&- \int_0^\tau \int_{\Omega_F} \mathbb{S} (\Ds \vu) : \Grad \bfphi  \dx \dt +
	\int_0^\tau \int_{\Omega_F} {\mu_F} \Big[ \vc{H} \otimes \vc{H} - \frac{1}{2} |\vc{H}|^2 \mathbb{I} \Big] : \Grad \bfphi  \dx \dt
	\label{Z21} 		
\end{align}	
for any $0 \leq \tau \leq T$ and any $\bfphi \in C^1_c(\Omega_F; R^d)$.

\end{Proposition}	

\section{Main results}
\label{M}

Having established the convergence in the fluid equations, we are ready to state our main results. 
As the {limits} of the equation of continuity as well as the momentum equation have been identified in Proposition \ref{ZP1}, {the next step is to} introduce the weak formulation of the limit problem for the induction equation.

\subsection{Weak form of the limit induction equation}
\subsubsection{Perfect electric isolator}
\label{pei}

The perfect electric isolator will be identified with the singular limit
\begin{equation} \label{B9} 
	\eta^\ep_{\rm ext} \to \infty.
\end{equation}
Recall that the boundary conditions on $\partial \Omega$ read
\begin{equation} \label{B10}
	\vc{H} \times \vc{n} = \vc{H}_{\rm ext} \times \vc{n},\ 
	\mu_F \vc{H} \cdot \vc{n} = \mu_{\rm ext} \vc{H}_{\rm ext} \cdot \vc{n}
	\ \mbox{in}\ (0,T) \times \partial \Omega,
\end{equation}	
where
\begin{align}
\Curl \vc{H}_{\rm ext} = 0 ,\ \Div (\mu_{\rm ext} \vc{H}_{\rm ext} ) = 0	\ \mbox{in}\ (0,T) \times \Omega_{\rm ext}.
\nonumber
\end{align}

The relevant weak formulation of the {limit} induction equation reads
\begin{equation} \label{B10a}
\Div (\mu \vc{H}) = 0,\ \mu(x) = \left\{ \begin{array}{l} \mu_{\rm int} \ \mbox{if}\ x \in \Omega_{\rm int}\\ 
\mu_F \ \mbox{if}\ x \in \Omega_F, \\ 
\mu_{\rm ext} \ \mbox{if}\ x \in {\Omega_{\rm ext}}	\end{array} \right.
\end{equation}
\begin{equation} \label{B10b}
\vc{H} \in L^\infty(0,T; L^2(\Td; R^d)), 
\ \Curl \vc{H} \in L^2(0,T; L^2(\Td; R^d)),\ \Curl \vc{H}|_{\Omega_{\rm ext}} = 0,
\end{equation}
and the integral identity
\begin{equation} \label{B10c}
	\intTd{ \mu \vc{H}(\tau, \cdot) \cdot \bfphi} = 
	\intTd{\mu \vc{H}_0 \cdot \bfphi }
	-
	\int_0^\tau \intO{ \Big[  
		(\mu \vc{H} \times \vu) \cdot \Curl \bfphi + \eta \Curl \vc{H} \cdot \Curl \bfphi \Big] } \dt 
\end{equation}		
for any $0 \leq \tau \leq T$, and any $\bfphi \in C^1(\Td)$, $\Curl \bfphi |_{\Omega_{\rm ext}} = 0$.

\subsubsection{Perfect magnetic conductor (PMC)}

The corresponding boundary conditions are of Dirichlet type
\begin{equation} \label{B12}
	\vc{H} \times \vc{n}|_{\partial \Omega} = 0.
\end{equation}	
The relevant weak formulation reads
\begin{equation} \label{B11c}
	\intO{ \mu \vc{H}(\tau, \cdot) \cdot \bfphi} = 
	\intO{\mu \vc{H}_0 \cdot \bfphi }
	-
	\int_0^\tau \intO{ \Big[  
		(\mu \vc{H} \times \vu) \cdot \Curl \bfphi + \eta \Curl \vc{H} \cdot \Curl \bfphi \Big] } \dt 
\end{equation}		
for any $0 \leq \tau \leq T$, and any $\bfphi \in C^1_c(\Omega; R^d)$.
{In addition, we have}
\[
{ \intO{ \mu \vc{H}(\tau, \cdot) \cdot \Grad \psi} = 0 \ 
	\mbox{whenever}\ \intO{ \mu \vc{H}_0 \cdot \Grad \psi} = 0}
\]
{meaning}
\begin{equation} \label{B11d}
{\Div (\mu \vc{H}) = 0}.	
\end{equation}

\subsubsection{Perfect electric conductor (PEC)}

The boundary conditions read
\begin{equation} \label{B13}
	\vc{H} \cdot \vc{n}|_{\partial \Omega} = 0,\ 
	\Curl \vc{H} \times \vc{n}|_{\partial \Omega} = 0.
\end{equation}	 
The weak formulation is 
\begin{equation} \label{B101a}
	\Div (\mu \vc{H}) = 0,\ \mu(x) = \left\{ \begin{array}{l} \mu_{\rm int} \ \mbox{if}\ x \in \Omega_{\rm int}\\ 
		\mu_F \ \mbox{if}\ x \in \Omega_F, 	\end{array} \right.
\end{equation}
and
\begin{equation} \label{B13a}
	\intO{ \mu \vc{H}(\tau, \cdot) \cdot \bfphi} = 
	\intO{\mu \vc{H}_0 \cdot \bfphi }
	-
	\int_0^\tau \intO{ \Big[  
		(\mu \vc{H} \times \vu) \cdot \Curl \bfphi + \eta \Curl \vc{H} \cdot \Curl \bfphi \Big] } \dt 
\end{equation}
for any $0 \leq \tau \leq T$, and any $\bfphi \in C^1(\Ov{\Omega}; R^d)$.

\subsubsection{Isolator type boundary conditions}

The boundary conditions read
\begin{equation} \label{B15}
\vc{H} \cdot \vc{n}|_{\partial \Omega} = 0,\ 
\Curl \vc{H} \cdot \vc{n}|_{\partial \Omega} = 0.	
	\end{equation}
The weak formulation is
\begin{equation} \label{B15b}
	\Div (\mu \vc{H}) = 0,\ \mu(x) = \left\{ \begin{array}{l} \mu_{\rm int} \ \mbox{if}\ x \in \Omega_{\rm int}\\ 
		\mu_F \ \mbox{if}\ x \in \Omega_F, 	\end{array} \right. ,
\end{equation}
and
 \begin{equation} \label{B15c}
 	\intO{ \mu \vc{H}(\tau, \cdot) \cdot \bfphi} = 
 	\intO{\mu \vc{H}_0 \cdot \bfphi }
 	-
 	\int_0^\tau \intO{ \Big[  
 		(\mu \vc{H} \times \vu) \cdot \Curl \bfphi + \eta \Curl \vc{H} \cdot \Curl \bfphi \Big] } \dt 
 \end{equation}		
 for any $0 \leq \tau \leq T$, and any $\bfphi \in C^1_c(\Omega; R^d)$.

\subsection{Singular limits - main results}

Having established the necessary preliminary material, we may state the main results of the present paper. 

\begin{Theorem}[\bf Singular limits] \label{MT1}
	Let $\Omega_{\rm int}$, $\Omega_F$, $\Omega_{\rm ext}$ satisfy the hypotheses {\eqref{i5}--\eqref{i6b}}. 	
	Let
	\begin{align} 
		\beta^\ep_F &= 0,\ \beta^\ep_{\rm int} \to \infty, \
		\beta^\ep_{\rm ext} \to \infty, \br 
		\nu^\ep_F &= \nu_F > 0,\ \lambda^\ep_F = \lambda_F \geq 0,\br 
		\nu^\ep_{\rm int} &\to \infty,\ \nu^\ep_{\rm ext} \to \infty \ \mbox{if}\ d = 3, \
		\lambda^\ep_{\rm int} \to \infty,\ \lambda^\ep_{\rm ext} \to \infty 
		\ \mbox{if}\ d = 2, \br 
		\mu^\ep_{\rm int} &= \mu_{\rm int} > 0, 
		\eta^\ep_{\rm int} = \eta_{\rm int} > 0,\ \mu^\ep_F = \mu_F > 0,\ 
		\eta^\ep_F = \eta_F > 0
		\label{BB1}
	\end{align}
{as $\ep \to 0$.}

	Then the following holds: 
	
	\begin{itemize}
		\item{\bf Perfect electric isolator.}
		In addition to \eqref{BB1}, suppose	
		\begin{equation} \label{BB2}
			\mu^\ep_{\rm ext} = \mu_{\rm ext} > 0,\  \eta^\ep_{\rm ext} \to \infty \ {\mbox{as}\ \ep \to 0}.
		\end{equation}
		Let $\vre$, $\vue$, $\vc{H}_\ep$ be a dissipative weak solution {of} the primitive compressible MHD system, with the initial data
		\begin{align} 
			\vr_{0,\ep} &= \vr_0,\ \int_{\Omega_F} \vr_0 \dx > 0, \ \vm_{0,\ep} = \vm_{0} \ \mbox{in}\ {\Td},\ \int_{\Omega_{\rm int} \cup \Omega_{\rm ext}} \frac{|\vm_{0,\ep}|^2}{\vr_0} \dx \to 0,\br
			\vc{H}_{0,\ep} &\to \vc{H}_0 \ \mbox{in}\ L^2(\Td; R^d),\ \intTd{ \mu^\ep \vc{H}_{0,\ep} \cdot \Grad \phi } = 0 \ \mbox{for any}\ \phi \in C^1(\Td),\ \br   
			\intTd{ E(\vr_0, \vm_{0,\ep}, \vc{H}_{0,\ep} )} &\leq c
			\nonumber
		\end{align}
		uniformly for $\ep \to 0$.
		
		Then, up to a suitable subsequence,
		\begin{align}
			\vre &\to \vr \ \mbox{in}\ C([0,T]; L^\gamma(\Td)) \cap L^1((0,T) \times \Td), \br
			\vue &\to \vu \ \mbox{weakly in}\ L^2(0,T; W^{1,2}(\Td; R^d)), \br
			\vc{H}_\ep &\to \vc{H} \ \mbox{weakly-(*) in}\ L^\infty(0,T; L^2(\Td; R^d)), 	
			\label{BB3}
		\end{align}
		where the limit $(\vr, \vu, \vc{H})$ is a finite energy weak solution of the compressible MHD with perfect electric isolator boundary conditions 
		{ \eqref{Z20}, \eqref{Z21}, \eqref{B10a}--\eqref{B10c}} { in $(0,T) \times \Td$}.
		
		\item {\bf Perfect magnetic conductor (PMC)}
		
		In addition to \eqref{BB1}, suppose	
		\begin{equation} \label{BB4}
			\mu^\ep_{\rm ext} \to \infty,\ \eta^\ep_{\rm ext} = \eta_{\rm ext} > 0 \ {\mbox{as}\ \ep \to 0}.
		\end{equation}
		{Let $\vre$, $\vue$, $\vc{H}_\ep$ be a dissipative weak solution {of} the primitive compressible MHD system, with the initial data}
		\begin{align} 
			\vr_{0,\ep} &= \vr_0,\ \int_{\Omega_F} \vr_0 \dx > 0, \ \vm_{0,\ep} = \vm_{0} \ \mbox{in}\ \Omega_F,\ \int_{\Omega_{\rm int} \cup \Omega_{\rm ext}} \frac{|\vm_{0,\ep}|^2}{\vr_0} \dx \to 0,\br
			\vc{H}_{0,\ep} &\to \vc{H}_0 \ \mbox{in}\ L^2(\Td; R^d),\ \intTd{ \mu^\ep \vc{H}_{0,\ep} \cdot \Grad \phi } = 0 \ \mbox{for any}\ \phi \in C^1(\Td),\ \br   
			\int_{\Omega_{\rm ext}} \mu^\ep_{\rm ext} |\vc{H}_{0, \ep}|^2 \dx &\to 0,\  
			\intTd{ E(\vr_0, \vm_{0,\ep}, \vc{H}_{0,\ep} )} \leq c
			\nonumber
		\end{align}
		uniformly for $\ep \to 0$.
		
		Then, up to a suitable subsequence
		\begin{align}
			\vre &\to \vr \ \mbox{in}\ C([0,T]; L^\gamma(\Td)) \cap L^1((0,T) \times \Td), \br
			\vue &\to \vu \ \mbox{weakly in}\ L^2(0,T; W^{1,2}_0(\Omega_F; R^d)), \br
			\vc{H}_\ep &\to \vc{H} \ \mbox{weakly-(*) in}\ L^\infty(0,T; L^2(Q; R^d)),\ \mbox{for any compact}\ Q \subset {\Omega},	
			\label{BB5}
		\end{align}
		where the limit $(\vr, \vu, \vc{H})$ is a finite energy weak solution to the compressible MHD with perfect magnetic conductor boundary conditions 
		\eqref{Z20}, \eqref{Z21},  {\eqref{B12}--\eqref{B11d}} in $(0,T) \times \Omega$.
		
		\item {\bf Perfect electric conductor (PEC)}
		
		In addition to \eqref{BB1}, suppose	
		\begin{equation} \label{BB6}
			\mu^\ep_{\rm ext} \to 0,\ \eta^\ep_{\rm ext} \to 0 
			\ { \mbox{as}\ \ep \to 0.}
		\end{equation}
			{Let $\vre$, $\vue$, $\vc{H}_\ep$ be a dissipative weak solution {of} the primitive compressible MHD system, with the initial data}
		\begin{align} 
			\vr_{0,\ep} &= \vr_0,\ \int_{\Omega_F} \vr_0 \dx > 0, \ \vm_{0,\ep} = \vm_{0} \ \mbox{in}\ \Omega_F,\ \int_{\Omega_{\rm int} \cup \Omega_{\rm ext}} \frac{|\vm_{0,\ep}|^2}{\vr_0} \dx \to 0,\br
			\vc{H}_{0,\ep} &\to \vc{H}_0 \ \mbox{in}\ L^2(\Td; R^d),\ \intTd{ \mu^\ep \vc{H}_{0,\ep} \cdot \Grad \phi } = 0 \ \mbox{for any}\ \phi \in C^1(\Td),\ \br   
			\int_{\Omega_{\rm ext}} \mu^\ep_{\rm ext} |\vc{H}_{0, \ep}|^2 \dx &\to 0,\  
			\intTd{ E(\vr_0, \vm_{0,\ep}, \vc{H}_{0,\ep} )} \leq c
			\nonumber
		\end{align}
		uniformly for $\ep \to 0$.
		
		Then, up to a suitable subsequence,
		\begin{align}
			\vre &\to \vr \ \mbox{in}\ C([0,T]; L^\gamma(\Td)) \cap L^1((0,T) \times \Td), \br
			\vue &\to \vu \ \mbox{weakly in}\ L^2(0,T; W^{1,2}_0(\Omega_F; R^d)), \br
			\vc{H}_\ep &\to \vc{H} \ \mbox{weakly-(*) in}\ L^\infty(0,T; L^2(Q; R^d)),\ \mbox{for any compact}\ Q \subset {\Omega},	
			\label{BB7}
		\end{align}
		where the limit $(\vr, \vu, \vc{H})$ is a finite energy weak solution to the compressible MHD with perfect electric conductor boundary conditions {
		\eqref{Z20}, \eqref{Z21}, {\eqref{B13}--\eqref{B13a}} in $(0,T) \times \Omega$}.
		
		\item {\bf Isolator type boundary conditions.}
		
		In addition to \eqref{BB1}, suppose	
		\begin{equation} \label{BB8}
			\mu^\ep_{\rm ext} \to 0,\ \eta^\ep_{\rm ext} \to \infty 
			\ {\mbox{as}\ \ep \to 0.} 
		\end{equation}
	{Let $\vre$, $\vue$, $\vc{H}_\ep$ be a dissipative weak solution {of} the primitive compressible MHD system, with the initial data}		
		\begin{align} 
			\vr_{0,\ep} &= \vr_0,\ \int_{\Omega_F} \vr_0 \dx > 0, \ \vm_{0,\ep} = \vm_{0} \ \mbox{in}\ \Omega_F,\ \int_{\Omega_{\rm int} \cup \Omega_{\rm ext}} \frac{|\vm_{0,\ep}|^2}{\vr_0} \dx \to 0,\br
			\vc{H}_{0,\ep} &\to \vc{H}_0 \ \mbox{in}\ L^2(\Td; R^d),\ \intTd{ \mu^\ep \vc{H}_{0,\ep} \cdot \Grad \phi } = 0 \ \mbox{for any}\ \phi \in C^1(\Td),\ \br   
			\int_{\Omega_{\rm ext}} \mu^\ep_{\rm ext} |\vc{H}_{0, \ep}|^2 \dx &\to 0,\  
			\intTd{ E(\vr_0, \vm_{0,\ep}, \vc{H}_{0,\ep} )} \leq c
			\nonumber
		\end{align}
		uniformly for $\ep \to 0$.
		
		Then, up to a suitable subsequence,
		\begin{align}
			\vre &\to \vr \ \mbox{in}\ C([0,T]; L^\gamma(\Td)) \cap L^1((0,T) \times \Td), \br
			\vue &\to \vu \ \mbox{weakly in}\ L^2(0,T; W^{1,2}_0(\Omega_F; R^d)), \br
			\vc{H}_\ep &\to \vc{H} \ \mbox{weakly-(*) in}\ L^\infty(0,T; L^2(Q; R^d)),\ \mbox{for any compact}\ Q \subset {\Omega},	
			\label{BB9}
		\end{align}
		where the limit $(\vr, \vu, \vc{H})$ is a finite energy weak solution to the compressible MHD with isolator type boundary conditions 
		\eqref{Z20}, \eqref{Z21}, \eqref{B15}, \eqref{B15b}, \eqref{B15c} in $(0,T) \times \Omega$.
	\end{itemize}	
\end{Theorem}

The rest of the paper is devoted to the proof of Theorem \ref{MT1}. Given the fact that the limit in the continuity and momentum equations has already been performed in Proposition 
\ref{ZP1}, we focus on the limit in the induction equation. The reader will have noticed that the limit system satisfies the energy inequality only in the integrated form \eqref{P12}. This is the price to pay when using a penalization method. A similar drawback has been noticed in \cite{FeNeSt}. As a matter of fact, we could have assumed the same form of the energy inequality also for the primitive system.

\section{Proof of the main results -- asymptotic limit in the induction equation}
\label{A}

We peform, case by case, the asymptotic limit in the induction equation as claimed in Theorem \ref{MT1}.

\subsection{Perfect electric isolator}

We have
\begin{equation} \label{A14}
	{\mu^\ep_{\rm ext}} = \mu_{\rm ext} > 0 \ \mbox{fixed}, \ 
	\eta^\ep_{\rm ext} \to \infty. 
\end{equation}
Consequently, using the bounds provided by the energy inequality \eqref{P12}, we get 
\begin{align} 
	\vc{H}_\ep &\to \vc{H} \ \mbox{in}\ C_{\rm weak}([0,T]; L^2(\Td; R^d)), \br
	\Curl \vc{H}_\ep &\to \Curl \vc{H}\ \mbox{weakly in}\ L^2((0,T) \times \Td; R^d), \br
	\Curl \vc{H}_\ep &\to 0 = \Curl \vc{H} \ \mbox{in}\ 
	L^2((0,T) \times \Omega_{\rm ext}; R^d).
	\label{A15}
\end{align}

The limit function satisfies a weak formulation of the induction equation 
\[
\int_0^T \intTd{ \Big[ \mu \vc{H} \cdot \partial_t \bfphi - 
	(\mu \vc{H} \times \mathds{1}_{\Omega_F} \vu) \cdot \Curl \bfphi - \eta \Curl \vc{H} \cdot \Curl \bfphi \Big] } \dt = - \intTd{ \mu \vc{H}_0 \cdot \bfphi }
\]		
for any $\bfphi \in C^1_c([0,T) \times \Td; R^d)$, $\Curl \bfphi|_{\Omega_{\rm ext}} = 0$. Moreover, 
\[
\Div(\mu \vc{H}) = 0, \ \mu = \left\{ \begin{array}{l} \mu_{\rm int} 
	\ \mbox{in}\ \Omega_{\rm int}, \\ \mu_F \ \mbox{in}\ \Omega_F, \\ 
	\mu_{\rm ext} \ \mbox{in}\ \Omega_{\rm ext}.
\end{array} \right.
\]	

The limit passage in the energy inequality \eqref{P12} is a routine matter based on weak lower semi--continuty of convex functionals and the hypotheses imposed on the initial data.

\subsection{Perfect magnetic conductor (PMC)}

Now, we consider the scaling
\begin{equation} \label{A18}
	\mu^\ep_{\rm ext} \to \infty,\ \eta^\ep_{\rm ext} = \eta_{\rm ext} > 0 \ \mbox{fixed}. 
\end{equation}

Similarly to the preceding part, we deduce 
\begin{align} 
	\vc{H}_\ep &\to \vc{H} \ \mbox{in}\ C_{\rm weak}([0,T]; L^2(\Omega; R^d)), \br
	\vc{H}|_{\Omega_{\rm ext}} &= 0, \br
	\Curl \vc{H}_\ep &\to \Curl \vc{H}\ \mbox{weakly in}\ L^2((0,T) \times \Td; R^d).
	\label{A19}
\end{align}

Now, the limit satisfies the induction equation in $(0,T) \times \Omega$, 
\[
\int_0^T \intO{ \Big[ \mu \vc{H} \cdot \partial_t \bfphi - 
	(\mu \vc{H} \times \mathds{1}_{\Omega_F} \vu) \cdot \Curl \bfphi - \eta \Curl \vc{H} \cdot \Curl \bfphi \Big] } \dt = - \intO{ \mu \vc{H}_0 \cdot \bfphi }
\]		
for any $\bfphi \in C^1_c([0,T) \times \Omega; R^d)$, together with 
\[
\Div(\mu \vc{H}) = 0, \ \mu = \left\{ \begin{array}{l} \mu_{\rm int} 
	\ \mbox{in}\ \Omega_{\rm int}, \\ \mu_F \ \mbox{in}\ \Omega_F .
\end{array} \right.
\]	

Finally, as $\vc{H}$ vanishes outside $\Omega$ and $\Curl \vc{H}$ exists in $\Td$, we deduce the desired boundary condition 
\[
\vc{H} \times \vc{n}|_{\partial \Omega} = 0
\]
using continuity of the tagent traces.

\subsection{Perfect electric conductor}

We examine the situation 
\begin{equation} \label{A23}
	\mu^\ep_{\rm ext} \to 0,\ \eta^\ep_{\rm ext} \to 0.
\end{equation}
First observe 
\begin{equation} \label{A24}
	\mu_\ep \vc{H}_\ep \to \mu \vc{H}, \ \mu_{\rm ext} = 0
\end{equation} 
As $\mu_F > 0$, the continuity of the normal trace yields 
\[
\vc{H} \cdot \vc{n}|_{\partial \Omega} = 0.	 
\]

Seeing that $\eta^\ep_{\rm ext} \to 0$, one can apply the same argument to 
the induction equation obtaining
\begin{equation} \label{A27}
	\int_0^T \intO{ \Big[ \mu \vc{H} \cdot \partial_t \bfphi - 
		(\mu \vc{H} \times \mathds{1}_{\Omega_F} \vu) \cdot \Curl \bfphi - \eta \Curl \vc{H} \cdot \Curl \bfphi \Big] } \dt = - \intO{ \mu \vc{H}_0 \cdot \bfphi }
\end{equation}	
for any $\bfphi \in C^1_c([0,T) \times \Ov{\Omega}; R^d)$. This yields the Neumann type boundary condition 
\begin{equation} \label{A28}
	\Curl \vc{H} \times \vc{n}|_{\partial \Omega} = 0
\end{equation}
implicitly satisfied by the choice of test functions in \eqref{A27}.

\subsection{Isolator type boundary conditions}

Finally, we consider the limit 
\begin{equation} \label{A29}
	\mu^\ep_{\rm ext} \to 0,\ \eta^\ep_{\rm ext} \to \infty.
\end{equation}	 
Exactly as in the preceding step, we deduce the condition 
\begin{equation} \label{A30}
	\vc{H} \cdot \vc{n}|_{\partial \Omega} = 0.
\end{equation}

In contrast with the previous step, however, the induction equation
\begin{equation} \label{A31}
	\int_0^T \intO{ \Big[ \mu \vc{H} \cdot \partial_t \bfphi - 
		(\mu \vc{H} \times \mathds{1}_{\Omega_F} \vu) \cdot \Curl \bfphi - \eta \Curl \vc{H} \cdot \Curl \bfphi \Big] } \dt = - \intO{ \mu \vc{H}_0 \cdot \bfphi }
\end{equation}	
is satisfied only for compactly supported test functions
$\bfphi \in C^1_c([0,T) \times \Omega; R^d)$. 
However, $\Curl H_\ep$ is a solenoidal function converging weakly to a limit that vanishes on $\Omega_{\rm ext}$. By continuity of the normal trace, we get 
\[
	\Curl \vc{H} \cdot \vc{n}|_{\partial \Omega} = 0.
\]

\def\cprime{$'$} \def\ocirc#1{\ifmmode\setbox0=\hbox{$#1$}\dimen0=\ht0
	\advance\dimen0 by1pt\rlap{\hbox to\wd0{\hss\raise\dimen0
			\hbox{\hskip.2em$\scriptscriptstyle\circ$}\hss}}#1\else {\accent"17 #1}\fi}


\end{document}